\documentclass[12pt]{article}

\usepackage[utf8]{inputenc}
\usepackage[T1]{fontenc}
\usepackage[english]{babel}
\usepackage{amsmath,amssymb,amsfonts}
\usepackage{geometry}

\title{
Artificial Intelligence and the Autonomization of Mathematics\\[0.3em]
{\small
\shortstack{
To what extent does mathematical activity require something beyond\\
the sophisticated exploration of formal structures?
}}
}

\author{Jaime Ripoll}

\date{}

\begin{document}


\maketitle

\begin{abstract}

This essay examines the relationship between artificial intelligence and the
historical evolution of modern Mathematics. Rather than interpreting
artificial intelligence as an external rupture imposed upon mathematical
activity, we argue that its growing effectiveness reveals a structural
tendency already present in the autonomization of Mathematics itself.

The paper analyzes how modern Mathematics progressively constructed formal
environments that became increasingly autonomous, stable, and internally
navigable, gradually reducing their dependence on concrete experience. In
contrast with Physics, which remains inevitably constrained by empirical
reality, Mathematics evolved toward forms of internal structural legitimacy
based primarily on formal coherence.

Within this context, the remarkable affinity between artificial intelligence
and contemporary mathematical practice appears less accidental than it might
initially seem. Systems capable of symbolic manipulation, structural
reorganization, and systematic exploration of formal consequences operate
precisely in the type of environment progressively stabilized by modern
Mathematics.

The essay also discusses possible limits of formal navigability. The creation
of genuinely new conceptual regimes, the emergence of new forms of
mathematical intelligibility, and the role of intuition and partially
pre-formal exploration may still involve dimensions not entirely reducible to
automated structural exploration. Husserl's reflections on mathematization
and the distancing of science from the \emph{Lebenswelt} provide a broader
philosophical framework for understanding this process.

Finally, we suggest that the contemporary debate surrounding artificial
intelligence may concern less a threat to Mathematics itself than a challenge
to the historical image we have constructed of the mathematician as the
privileged interpreter of mathematical structures.

\end{abstract}

\section{Introduction}

The recent development of artificial intelligence has once again brought to
the forefront an old question: what exactly characterizes human mathematical
activity? At first sight, the issue appears to be mainly technological. Much
has been said about the automation of proofs, pattern recognition, automatic
generation of conjectures, and computational exploration of large formal
spaces. Nevertheless, the problem may be deeper than that.

Artificial intelligence seems to make particularly visible a tension that has
silently accompanied an important part of the modern history of Mathematics:
to what extent does mathematical activity effectively depend on human
subjectivity?

For a long time, questions of this kind would have appeared almost absurd.
Mathematics was traditionally regarded as one of the highest expressions of
human intelligence. Intuition, creativity, elegance, and conceptual depth
seemed inseparable from the figure of the mathematician.

Yet the very evolution of modern Mathematics may suggest a more complex
movement. An important part of this evolution consisted in the progressive
construction of formal universes that became increasingly autonomous, stable,
and internally navigable.

Throughout this process, Mathematics gradually came to operate on structures
in which internal relations, limits, and infinite processes exhibit highly
regularized behavior. This made it possible to explore large formal spaces
with increasing degrees of structural stability.

From this perspective, artificial intelligence may not represent an external
rupture with respect to Mathematics. Rather, it may represent the extreme
manifestation of a structural tendency that Mathematics itself historically
helped to construct.

\section{Structural Autonomization and Completeness}

In order to understand why artificial intelligence now finds a particularly
favorable environment in Mathematics, it may be necessary to look at a
broader historical process: the progressive autonomization of mathematical
structures with respect to concrete experience.

In the natural sciences, confrontation with empirical reality remains
unavoidable. A physical theory does not survive merely through logical
coherence. It must accommodate measurements, face experimental noise, and
withstand continuous confrontation with the physical world. Mathematically
elegant structures may simply turn out to be physically wrong.

Modern Mathematics developed progressively in a different direction. During a
large part of its history, the legitimacy of mathematical objects still
seemed to depend on some form of geometrical, intuitive, or physical support.
Over time, however, increasingly abstract structures gradually came to be
legitimized through internal formal coherence.

The history of Mathematics itself is full of examples of this kind: negative
numbers, complex numbers, infinitesimals, non-Euclidean geometries, and
infinite-dimensional spaces spent long periods being regarded almost as
intellectual fictions.

The initial reaction to non-Euclidean geometries is perhaps particularly
revealing. The problem did not consist merely in abandoning the parallel
postulate. The deeper difficulty was to accept a geometry that seemed to
break with the ordinary intuition of physical space.

Once axioms, definitions, and inferential rules were established, entire
theories began to acquire an internal dynamics relatively independent of
concrete experience. Structures initially introduced in order to describe
physical phenomena frequently revealed unexpected connections, autonomous
generalizations, and internal properties that far exceeded their original
motivation.

This autonomization did not produce merely greater abstraction. It also
contributed to the construction of formal environments endowed with strong
structural stability.

Some fundamental assumptions of modern Analysis make this particularly
visible. Consider, for example, a sequence of closed cubes, each contained in
the preceding one, whose diameters tend to zero. Modern Mathematics assumes
that the intersection of all these cubes reduces exactly to a single point.

Today this assumption appears almost natural. Yet it constitutes an extremely
strong idealization. Physically, it is not even clear what it would mean to
continue such a process indefinitely at arbitrarily small scales. Material
limitations, quantum effects, and restrictions of measurement render the
concrete situation far less regular than the idealized mathematical
continuum.

The point, however, is not to question the legitimacy of such constructions.
The relevant issue here is structural. Modern assumptions of completeness
produce extraordinarily stable formal universes in which infinite processes
admit guaranteed closure. Limits, continuity, compactness, spectral theory,
functional spaces, and a large part of modern mathematical Physics depend
precisely on this stability.

In a certain sense, modern Mathematics progressively constructed environments
that became increasingly coherent, complete, and internally navigable. And
environments of this kind favor precisely the systematic exploration of
structural relations.

\section{Artificial Intelligence and Formal Exploration}

The affinity between artificial intelligence and Mathematics becomes clearer
when we observe the nature of these progressively autonomized formal
environments.

Artificial systems operate with increasing efficiency precisely in tasks such
as symbolic manipulation, structural reorganization, pattern recognition, and
systematic exploration of formal consequences. And these capacities coincide
in a striking way with important parts of contemporary mathematical practice.

This does not necessarily imply understanding in the traditional human sense.
But it progressively makes less evident what exactly separates certain
sophisticated forms of mathematical activity from highly developed processes
of structural exploration.

The singularity of Mathematics makes this question particularly delicate. A
large part of mathematical activity consists precisely in exploring internal
relations within highly stabilized formal universes.

At this point, an uncomfortable distinction emerges between Mathematics and
the mathematician.

Mathematics seems to possess a structural objectivity relatively independent
of its historical explorers. Theorems remain true regardless of who proves
them. Mathematical relations persist independently of culture, biology, or
even the existence of particular human subjects.

The mathematician, by contrast, is contingent. His or her activity depends on
education, language, memory, cognitive limitations, and historical context.

The more mathematical objects came to be legitimized through internal
structural coherence, the less mathematical activity seemed to depend
directly on concrete experience. Artificial intelligence now makes visible a
particularly radical possibility: perhaps a significant part of mathematical
exploration depends less on human subjectivity than has traditionally been
assumed.

There is also an important difference between Mathematics and Physics in this
context. The greatest difficulty for artificial systems may not lie in the
exploration of already formalized mathematical structures, but rather in the
extraction of relevant structures from the empirical world. Physics remains
inevitably tied to confrontation with concrete reality. Mathematics, once
autonomized, appears capable of continuing to evolve in a much more internal
manner.

\section{The Limits of Formal Navigability}

Nevertheless, it would be premature to conclude that all mathematical
activity can be reduced to the efficient exploration of formalized
structures.

There exists in the history of Mathematics a component that is far more
difficult to describe formally: the creation of new conceptual regimes. The
formulation of entire research programs --- as in Riemann, Grothendieck, or
the emergence of category theory --- seems to involve something different
from the mere navigation of already stabilized formal spaces.

To explore a mathematical environment presupposes that its relevant objects,
fundamental rules, and criteria of coherence are already minimally
constituted. Creating new conceptual horizons is another operation entirely:
it involves transforming the very space within which future explorations will
become possible.

Moreover, human mathematical activity does not occur only within stabilized
formal environments. It involves error, hesitation, vague analogies,
imperfect spatial intuitions, and long periods of partially disordered
exploration. Many discoveries seem to emerge precisely from cognitive regions
in which clear structural navigability does not yet exist.

There may lie here an important limit to formal autonomization itself.

Edmund Husserl's reflections help illuminate this problem from a broader
perspective. In \emph{Die Krisis der europ\"aischen Wissenschaften und die
transzendentale Ph\"anomenologie}, Husserl observed that the progressive
mathematization of the world produced a growing separation between the formal
systems of science and the immediate horizon of lived experience
(\emph{Lebenswelt}).

Seen in this way, the autonomization of Mathematics does not constitute
merely an internal technical phenomenon of the discipline. It also
participates in a broader cultural movement of symbolic sedimentation and
progressive distancing from the lived world.

Artificial intelligence may represent a particularly extreme moment in this
process. Artificial systems operate precisely upon formalized representations
and structurally stabilized environments. Yet the decisive question remains
open: to what extent can the production of new horizons of mathematical
intelligibility be entirely absorbed by this type of formal navigability?

\section{Conclusion}

Artificial intelligence may not be introducing an external rupture into
Mathematics. It may merely be making particularly visible a structural
tendency that silently accompanied an important part of its modern evolution.

The more Mathematics emancipated itself from concrete experience in order to
acquire formal autonomy, the more it seemed to approach environments in which
structural exploration can occur partially independently of human
subjectivity.

None of this implies, however, that all mathematical activity can be reduced
to the automated exploration of stabilized structures. The formulation of new
conceptual regimes, the creation of new forms of intelligibility, and the
profound reorganization of the mathematical horizon itself may still involve
dimensions that escape pure formal navigability.

For a long time, mathematical activity was romanticized as an almost
mysterious manifestation of an essentially irreducible human creativity.
Intuition, conceptual depth, and genius were frequently treated as signs of a
singularity practically unattainable by any form of automated process.

Yet the very evolution of modern Mathematics may have silently prepared a
more ambiguous scenario. In a profound irony, the more the discipline
constructed autonomous, coherent, and internally navigable formal
environments, the more an important part of mathematical activity came to
approximate systematic structural exploration.

Perhaps artificial intelligence does not truly threaten Mathematics.
Perhaps it threatens only the historical image we constructed of ourselves as
its privileged interpreters.

\section*{References}

Edmund Husserl,
\emph{Die Krisis der europ\"aischen Wissenschaften und die transzendentale
Ph\"anomenologie: Eine Einleitung in die ph\"anomenologische Philosophie},
edited by Walter Biemel,
Martinus Nijhoff,
The Hague,
1954.

English translation:

Edmund Husserl,
\emph{The Crisis of European Sciences and Transcendental Phenomenology:
An Introduction to Phenomenological Philosophy},
translated by David Carr,
Northwestern University Press,
Evanston,
1970.

\bigskip

Jaime Ripoll

Federal University of Rio Grande do Sul

Porto Alegre, Brazil

jaime.ripoll@ufrgs.br

\end{document}